\documentclass[a4paper,12pt]{article}
\usepackage{amssymb}

\usepackage{amsmath}

\input{tcilatex}

\begin{document}

\title{Multiplicative rule in the Grothendieck cohomology of a flag variety\\
{\small Dedicated to Professor Yifeng Sun on his 80th birthday }}
\author{Haibao Duan \\
Institute of Mathematics, Chinese Academy of Sciences\\
Beijing 100080, dhb@math.ac}
\date{}
\maketitle

\begin{abstract}
In the Grothendieck cohomology of a flag variety $G/H$ there are two
canonical additive bases, namely, the Demazure basis [D] and the
Grothendieck basis [KK].

We present explicit formulae that reduce the multiplication of these basis
elements to the Cartan numbers of $G$.
\end{abstract}

\section{Introduction}

Let $G$ be a compact connected Lie group and let $H\subset G$ be the
centralizer of a one--parameter subgroup in $G$. The homogeneous space $G/H$
is known as a \textsl{flag variety}. We fix a maximal torus $T\subseteq H $
and write $W$ (resp. $W^{^{\prime }}$) for the Weyl group of $G$ (resp. $H$).

In the founding article [AH] of topological $K$--theory as a ``generalized
cohomology theory'', Atiyah and Hirzebruch raised also the problem to
determine the ring $K(G/H)$, with the expectation that ``\textsl{the new
cohomology theory can be applied to various topological questions and may
give better results than the ordinary cohomology theory}''. As an initial
step they showed that the $K(G/H)$ is a free $\mathbb{Z}$-module with rank
equal to the quotient of the order of $W$ by the order of $W^{^{\prime }}$
[AH, Theorem 3.6].

In the subsequent years two distinguished additive bases of $K(G/H)$ have
emerged from algebraic geometry\footnote{%
The Grothendieck cohomology of coherent sheaves on $G/H$ is canonically
isomorphic to the $K$-theory of complex vector bundles on $G/H$ [Br$_{2}$,PR$%
_{2}$].}. The first of these, valid for the case $H=T$ and indexed by
elements from the Weyl group: $\{a_{w}\in K(G/T)\mid w\in W\}$, is due to
Demazure [D, Proposition 7]. The basis element $a_{w}$ will be called the 
\textsl{Demazure class} relative to $w\in W$.

The second basis goes back to the classical works of Grothendieck and
Chevalley. Let $l:W\rightarrow \mathbb{Z}$ is the length function relative
to a fixed maximal torus $T\subset G$ and identify the set $W/W^{\prime }$
of left cosets of $W^{\prime }$ in $W$ with the subset $\overline{W}=\{w\in
W\mid l(w)\leq l(w_{1})$ for all $w_{1}\in wW^{\prime }\}$ of $W$. According
to Chevalley [Ch] the flag variety $G/H$ admits a canonical partition into
Schubert varieties, indexed by elements of $\overline{W}$,

\begin{center}
$G/H=\underset{w\in \overline{W}}{\cup }X_{w}(H)$, $\quad \dim
X_{w}(H)=2l(w) $.
\end{center}

\noindent The coherent sheaves $\Omega _{w}(H)\in K(G/H)$ of the Schubert
variety Poincare dual to the $X_{w}(H)$, $w\in \overline{W}$, form a basis
for the $\mathbb{Z}$-module $K(G/H)$ by Grothendieck (cf. [C, Lecture 4]).
The $\Omega _{w}(H)$ is called the \textsl{Grothendieck class} relative to $%
w\in \overline{W}$.

\bigskip

A complete description of the ring $K(G/T)$ (resp. $K(G/H)$) now amounts to
specify the structure constants $C_{u,v}^{w}\in \mathbb{Z}$ (resp. $%
K_{u,v}^{w}(H)\in \mathbb{Z}$) required to express the products of basis
elements

\begin{enumerate}
\item[(1.1)] $\qquad a_{u}\cdot a_{v}=\sum C_{u,v}^{w}a_{w}$ (resp. $\Omega
_{u}(H)\cdot \Omega _{v}(H)=\sum K_{u,v}^{w}(H)\Omega _{w}(H)$),
\end{enumerate}

\noindent where $u,v,w\in W$ (resp. $\in \overline{W}$).

Based on combinatorial methods, partial information concerning the ring $%
K(G/H)$ has been achieved during the past decade. Generalizing the classical
Pieri-Chevalley formula in the ordinary cohomology [Ch], various
combinatorial formulae expressing the product $L\cdot \Omega _{u}$ in terms
of the $\Omega _{w}$ were obtained, where $L$ is a line bundle on $G/T$.
This was originated by Fulton and Lascoux [FL] for the unitary group $U(n)$
of rank $n$ (see also Lenart [L]), extended to general $G$ by Pittie and Ram
[PR$_{1,2}$], Mathieu [M], Littelmann and Seshadri [LS] by using very
different methods. Another progress is when $G=U(n)$ and $H=U(k)\times
U(n-k) $, the flag variety $G/H$ is the Grassmannian $G_{n,k}$ of $k$-planes
through the origin in $\mathbb{C}^{n}$ and a combinatorial description for
the $K_{u,v}^{w}(H)$ was obtained by Buch [Bu].

In general, using purely geometric approach, Brion [Br$_{1}$,Br$_{2}$]
proved that the $K_{u,v}^{w}$ have alternating signs, confirming a
conjecture of Buch in [Bu].

\bigskip

In this paper we present a formula that expresses the constants $C_{u,v}^{w}$
(resp. the $K_{u,v}^{w}(H)$ ) by the Cartan numbers of $G$, see Theorem 1 in
\S 2 (resp. Theorem 2 in \S 5). These results are natural generalization of
the formula in [Du$_{2}$] for multiplying Schubert classes in the cohomology
of $G/H$ for, in the special cases $l(w)=l(u)+l(v)$, the number $C_{u,v}^{w}$
(resp. the $K_{u,v}^{w}(H)$ ) agrees with the coefficient of the Schubert
class $P_{w}$ in the product $P_{u}\cdot P_{v}$ (see in the notes after
Theorem 1 in \S 2).

An important problem in algebraic combinatorics is to find a combinatorial
description of the $C_{u,v}^{w}$ (resp. the $K_{u,v}^{w}(H)$ ) that has the
advantage to reveal their signs [L,Bu]. Apart from the combinatorial concern,%
\textsl{\ }effective computation in the $K$--theory of such classical
manifolds as the $G/H$\ is decidedly required by many problems from geometry
and topology. We emphasis that the formulae as given in this paper are
computable, although it is not readily to reveal the signs of the constants.
With some additional works the existing program [DZ$_{2}$] for multiplying
Schubert classes can be extended to implement the $C_{u,v}^{w}$ (resp. $%
K_{u,v}^{w}(H)$). This program uses the Cartan matrix of $G$ as the only
input and, therefore, is functional uniformly for all $G/H$.

\bigskip

The author feels very grateful to S. Kumar for valuable communication.
Indeed, our exposition benefits a lot from certain results (cf. Lemma 4.2;
Lemma 5.1) developed in the classical treatise [KK] on this subject.

Thanks are also due to my referee for many improvements on the earlier
version of the paper, and for his kindness in informing me the work [GR] by
S. Griffeth and A. Ram, where a combinatorial method to multiply two
elements of the Grothendieck basis of the equivariant K--ring of a flag
variety is given, and the article [W$_{2}$] by M. Willems, where similar
methods are used in the setting of equivariant K--theory.

\section{The formula for $C_{u,v}^{w}$}

We introduce notations (from Definition 1 to 4) in terms of which the
formula for $C_{u,v}^{w}$ will be presented in Theorem 1.

\bigskip

Fix once and for all a maximal torus $T$ in $G$. Set $n=\dim T$. Equip the
Lie algebra $L(G)$ with an inner product $(,)$ so that the adjoint
representation acts as isometries of $L(G)$.

The restriction of the exponential map $\exp :L(G)\rightarrow G$ to $L(T)$
defines a set $D(G)$ of $\frac{1}{2}(\dim G-n)$ hyperplanes in $L(T)$, i.e.
the set of\textsl{\ singular hyperplanes\ }through the origin in $L(T)$. The
reflections of $L(T)$ in these planes generate the Weyl group $W$ of $G$
([Hu, p.49]).

Take a regular point $\alpha \in L(T)$ and let $\Delta $ be the set of
simple roots relative to $\alpha $ [Hu, p.47]. For a $\beta \in \Delta $ the
reflection $r_{\beta }$ in the hyperplane $L_{\beta }\in D(G)$ relative to $%
\beta $ is called a \textsl{simple reflection}. If $\beta ,\beta ^{\prime
}\in \Delta $, \textsl{the Cartan number\ }

\begin{center}
$\beta \circ \beta ^{\prime }=2(\beta ,\beta ^{\prime })/(\beta ^{\prime
},\beta ^{\prime })$
\end{center}

\noindent is always an integer (only $0,\pm 1,\pm 2,\pm 3$ can occur). It is
also customary to use $(\beta ,\beta ^{\prime \vee })$ instead of $\beta
\circ \beta ^{\prime }$.

\bigskip

It is well known that the set of simple reflections $\{r_{\beta }\mid \beta
\in \Delta \}$ generates $W$. That is, any $w\in W$ admits a factorization
of the form

\begin{enumerate}
\item[(2.1)] \ \ \ \ \ \ \noindent $w=r_{\beta _{1}}\cdot \cdots \cdot
r_{\beta _{m}}$,$\quad \beta _{i}\in \Delta $.
\end{enumerate}

\begin{quote}
\textbf{Definition 1.}\textit{\ }The \textsl{length} $l(w)$\ of a $w\in W$\
is the least number of factors in all decompositions of $w$\ in the form
(2.1). The decomposition (2.1) is said \textsl{reduced} if $m=l(w)$.

If (2.1) is a reduced decomposition, the $m\times m$\ (strictly upper
triangular) matrix $A_{w}=(a_{i,j})$\ with
\end{quote}

\begin{center}
$a_{i,j}=\{%
\begin{array}{c}
0\text{ \ if }i\geq j\text{;\qquad } \\ 
\beta _{j}\circ \beta _{i}\text{\ if }i<j%
\end{array}%
$
\end{center}

\begin{quote}
\noindent will be called \textsl{the Cartan matrix of }$w$\ relative to the
decomposition (2.1).

\textbf{Definition 2.} Given a sequence $(\beta _{1},\cdots ,\beta _{m})$\
of simple roots and a $w\in W$, let $[i_{1},\cdots ,i_{k}]\subseteq \lbrack
1,\cdots ,m]$\ be the subsequence maximal in the inverse-lexicographic order
so that
\end{quote}

\begin{center}
$l(w)>l(wr_{\beta _{i_{k}}})>$\ $l(wr_{\beta _{i_{k}}}r_{\beta
_{i_{k-1}}})>\cdots >l(wr_{\beta _{i_{k}}}r_{\beta _{i_{k-1}}}\cdots
r_{\beta _{i_{1}}})$.
\end{center}

\begin{quote}
\noindent We call $(\beta _{1},\cdots ,\beta _{m})$\ \textsl{a derived
(simple root) sequence of }$w$, written $(\beta _{1},\cdots ,$\ $\beta
_{m})\thicksim w$, if $k=l(w)$\ (i.e. $r_{\beta _{i_{1}}}\cdots r_{\beta
_{i_{k}}}$\ is a reduced decomposition of $w$).
\end{quote}

\textbf{Remark 1.} It is clear that $(\beta _{1},\cdots ,\beta
_{m})\thicksim w$ implies $m\geq l(w)$, while the equality holds if and only
if $w=r_{\beta _{1}}\cdots r_{\beta _{m}}$.

The Definition 2 implies also that if $e\in W$ is the group unit, then $%
(\beta _{1},\cdots ,\beta _{m})\sim e$ for any sequence $(\beta _{1},\cdots
,\beta _{m})$\ of simple roots.

\bigskip

Let $\mathbb{Z}[y_{1},\cdots ,y_{m}]$ be the ring of integral polynomials in 
$y_{1},\cdots ,y_{m}$, graded by $\mid y_{i}\mid =1$, and let $\mathbb{Z}%
[y_{1},\cdots ,y_{m}]_{(n)}$ be the submodule of all polynomials of degree $%
\leq n$. We introduce the additive maps

\begin{center}
$_{(n)}:\mathbb{Z}[y_{1},\cdots ,y_{m}]\rightarrow \mathbb{Z}[y_{1},\cdots
,y_{m}]_{(n)}$, $n\geq 0$,
\end{center}

\noindent by the following rule. If $f\in \mathbb{Z}[y_{1},\cdots ,y_{m}]$,
then

\begin{center}
$f=f_{(n)}+$ a sum of monomials of degree $>n$.
\end{center}

\begin{quote}
\textbf{Definition 3.} Let $A=(a_{i,j})_{m\times m}$\ be a strictly upper
triangular matrix (with integer entries). In terms of the entries of $A$\
define two sequences $\{q_{1},\cdots ,q_{m}\},$\ $\{\overline{q}_{1},\cdots ,%
\overline{q}_{m}\}\subset \mathbb{Z}[y_{1},\cdots ,y_{m}]$\ of polynomials
inductively as follows. Put $q_{1}=\overline{q}_{1}=1$,\ and for $k>1$\ let
\end{quote}

\begin{center}
$\qquad
q_{k}=\prod\limits_{a_{i,k}>0}(y_{i}+1)^{a_{i,k}}\prod%
\limits_{a_{i,k}<0}(-q_{i}y_{i}+1)^{-a_{i,k}}$;

$\qquad \overline{q}_{k}=\underset{a_{i,k}>0}{\prod }%
(-q_{i}y_{i}+1)^{a_{i,k}}\underset{a_{i,k}<0}{\prod }(y_{i}+1)^{-a_{i,k}}$.
\end{center}

\textbf{Remark 2.} As an example if $A=\left( 
\begin{array}{ccc}
0 & 1 & 2 \\ 
0 & 0 & -1 \\ 
0 & 0 & 0%
\end{array}%
\right) $, then

\begin{center}
$\qquad \{%
\begin{array}{c}
q_{2}=y_{1}+1,\  \\ 
\overline{q_{2}}=1-y_{1}\text{;}%
\end{array}%
$ and $\{%
\begin{array}{c}
q_{3}=(y_{1}+1)^{2}[-(y_{1}+1)y_{2}+1], \\ 
\overline{q_{3}}=[-(y_{1}+1)y_{1}+1]^{2}(y_{2}+1)\text{.}%
\end{array}%
$
\end{center}

\noindent Note also that, since $A$ is strictly upper triangular, we always
have

\begin{center}
$q_{k},\overline{q}_{k}\in \mathbb{Z}[y_{1},\cdots ,y_{k-1},q_{1},\cdots
,q_{k-1}]=\mathbb{Z}[y_{1},\cdots ,y_{k-1}]$.
\end{center}

\begin{quote}
\textbf{Definition 4.} Given a strictly upper triangular matrix $%
A=(a_{i,j})_{m\times m}$\ of rank $m$\ define the operator $\Delta _{A}:%
\mathbb{Z}[y_{1},\cdots ,y_{m}]_{(m)}\rightarrow \mathbb{Z}$\ as the
composition
\end{quote}

\begin{center}
$\mathbb{Z}[y_{1},\cdots ,y_{m}]_{(m)}\overset{D_{m-1}}{\rightarrow }\mathbb{%
Z}[y_{1},\cdots ,y_{m-1}]_{(m-1)}\overset{D_{m-2}}{\rightarrow }\cdots 
\overset{D_{1}}{\rightarrow }\mathbb{Z}[y_{1}]_{(1)}\overset{D_{0}}{%
\rightarrow }\mathbb{Z}$,
\end{center}

\begin{quote}
\noindent where the operator $D_{k-1}:\mathbb{Z}[y_{1},\cdots
,y_{k}]_{(k)}\rightarrow \mathbb{Z}[y_{1},\cdots ,y_{k-1}]_{(k-1)}$\ is
given by the following elimination rule.

\textsl{Expand each }$f\in \mathbb{Z}[y_{1},\cdots ,y_{k}]_{(k)}$\textsl{\
in terms of the powers of }$y_{k}$
\end{quote}

\begin{center}
$f=h_{0}+h_{1}y_{k}+h_{2}y_{k}^{2}+\cdots +h_{k}y_{k}^{k},$\textsl{\ \ }$%
h_{i}\in \mathbb{Z}[y_{1},\cdots ,y_{k-1}],$
\end{center}

\begin{quote}
\textsl{\noindent then put}
\end{quote}

\begin{center}
$D_{k-1}(f)=[h_{1}+h_{2}(\overline{q_{k}}-1)+\cdots +h_{k}(\overline{q_{k}}%
-1)^{k-1}]_{(k-1)}$\textsl{,}
\end{center}

\begin{quote}
\textsl{\noindent where the }$\overline{q_{k}}$\textsl{\ is given by }$A$%
\textsl{\ as in Definition 3 (see also Remark 2).}
\end{quote}

\textbf{Remark 3. }The $\Delta _{A}$ can be easily evaluated, as the formula
shows.

\begin{center}
$D_{k-1}(f)=[\sum\limits_{n\geq 1}(\frac{1}{n!}\frac{\partial ^{n}f}{%
(\partial y_{k})^{n}}\mid _{y_{k}=0})(\overline{q_{k}}-1)^{n-1}]_{(k-1)}$.
\end{center}

\textbf{Remark 4.} The operator\textbf{\ }$\Delta _{A}$ extends the idea of 
\textsl{triangular operator} $T_{A}$ in [Du$_{1}$,Du$_{2}$,DZ$_{1}$,DZ$_{2}$%
] in the following sense. If $f\in \mathbb{Z}[y_{1},\cdots ,y_{m}]$\ is of
homogeneous degree $m$, then $\Delta _{A}(f)=T_{-A}(f)$\textit{.}

The operator $T_{A}$ appears to be a useful tool in computing with the
cohomology of $G/H$. It was applied to evaluate the degrees of Schubert
varieties in [Du$_{1}$], to present a formula for multiplying Schubert
classes in [Du$_{2}$,DZ$_{2}$], and to compute the Steenrod operations on
the $\mathbb{Z}_{p}$--cohomologies of $G/H$ \ in [DZ$_{1}$]. Apart from the $%
\Delta _{A}$, another generalization of the $T_{A}$ was given by Willems [W,
Definition 5.2.1], which was useful for multiplying Schubert classes in the $%
T$-equivariant cohomology of $G/T$ [W$_{1}$, Theorem 5.3.1].

\bigskip

Assume that $w=r_{\beta _{1}}\cdots r_{\beta _{m}}$, $\beta _{i}\in \Delta $%
, is a reduced decomposition of $w\in W$, and let $A_{w}=(a_{i,j})_{m\times
m}$ be the associated Cartan matrix. For a subsequence $L=[i_{1},\cdots
,i_{k}]\subseteq \lbrack 1,\cdots ,m]$ we set

\begin{center}
$\beta (L)=(\beta _{i_{1}},\cdots ,\beta _{i_{k}})$,\qquad\ $%
y_{L}=y_{i_{1}}\cdots y_{i_{k}}\in \mathbb{Z}[y_{1},\cdots ,y_{m}]$.
\end{center}

\begin{quote}
\textbf{Theorem 1.} \textsl{For }$u,v\in W$\textsl{\ we have}
\end{quote}

\begin{enumerate}
\item[(2.2)] $\qquad (-1)^{l(w)}C_{u,v}^{w}=(-1)^{l(u)+l(v)}\Delta
_{A_{w}}[(\sum\limits_{\beta (L)\thicksim u}y_{L})(\sum\limits_{\beta
(K)\thicksim v}y_{K})]_{(m)}$

$\qquad \qquad \qquad \qquad \quad -\sum\limits_{\substack{ l(u)+l(v)\leq
l(x)\leq l(w)-1  \\ \beta (1,\cdots ,m)\thicksim x\in W}}%
(-1)^{l(x)}C_{u,v}^{x}$\textsl{,}
\end{enumerate}

\begin{quote}
\noindent \textsl{where }$L,K\subseteq \lbrack 1,\cdots ,m]$.
\end{quote}

As suggested by Theorem 1, the job to compute all $C_{u,v}^{x}$ for given $%
u,v\in W$ may be organized as follows.

\begin{enumerate}
\item[(1)] If $l(w)(=m)<l(u)+l(v)$, then

$\qquad \qquad \lbrack (\sum\limits_{\beta (L)\thicksim
u}y_{L})(\sum\limits_{\beta (K)\thicksim v}y_{K})]_{(m)}=0$

implies that $C_{u,v}^{w}=0$.

\item[(2)] If $l(w)(=m)=l(u)+l(v)$ the formula becomes
\end{enumerate}

\ \ \ \ $\qquad C_{u,v}^{w}=\Delta _{A_{w}}[(\sum\limits_{\beta (L)\thicksim
u}y_{L})(\sum\limits_{\beta (K)\thicksim v}y_{K})]_{(m)}$

$\qquad =\Delta _{A_{w}}[(\sum\limits_{r_{L}=u,\mid L\mid
=l(u)}y_{L})(\sum\limits_{r_{K}=v,\mid K\mid =l(v)}y_{K})]$ (cf. Remark 1)

$\qquad =T_{-A_{w}}[(\sum\limits_{r_{L}=u,\mid L\mid
=l(u)}y_{L})(\sum\limits_{r_{K}=v,\mid K\mid =l(v)}y_{K})]$ (cf. Remark 4),

\begin{enumerate}
\item[ ] \noindent where $r_{L}=r_{\beta _{i_{1}}}\cdots r_{\beta _{i_{k}}}$%
, $\mid L\mid =k$ if $L=[i_{1},\cdots ,i_{k}]$. This recovers the formula
for multiplying Schubert classes in the ordinary cohomology of $G/T$ [Du$%
_{2} $,DZ$_{1}$].

\item[(3)] In general, assuming that all the constants $C_{u,v}^{x}$ with $%
l(x)<m$ have been obtained, the theorem gives $C_{u,v}^{w}$ with $l(w)=m$ in
terms of the operator $\Delta _{A_{w}}$ as well as those $C_{u,v}^{x}$ ($%
l(x)<m$, $(\beta _{1},\cdots ,\beta _{m})\thicksim x$) calculated before.
\end{enumerate}

\noindent It is clear from the discussion that Theorem 1 reduces the $%
C_{u,v}^{w}$ to the operators $\Delta _{A_{x}}$, hence to the matrices $%
A_{x} $ formed by Cartan numbers.

\bigskip

Theorem 1 is originated from the celebrated Bott-Samelson cycles on flag
manifolds [BS]. This may be seen from the geometric consideration that
underlies the algebraic formation from Definition 1 to 4. Indeed, the Cartan
matrix of a $w$ (Definition 1) with respect to the decomposition (2.1) gives
the structural data characterizing the Bott-Samelson cycle $S(\alpha ;\beta
_{1},\cdots ,\beta _{m})$ as a twisted products of $2$--spheres (Lemma 4.3);
the polynomials $\overline{q}_{k}$'s ( Definition 3) provide the relations
in describing the $K$--ring of $S(\alpha ;\beta _{1},\cdots ,\beta _{m})$ as
the quotient of a polynomial ring (Lemma 4.4); the operator $\Delta _{A}$
(Definition 4) handles the integration along the top cell of $S(\alpha
;\beta _{1},\cdots ,\beta _{m})$ in the $K$--theory (Lemma 4.4); and the
idea of derived sequence of a Weyl group element(Definition 2) is required
to specify the induced map of a Bott-Samelson cycle in $K$--theory (Lemma
4.5).

The remaining sections are so arranged. Before involving the specialities of
flag manifolds, Section 3 studies the $K$-theory of \textsl{twisted products
of }$2$\textsl{-spheres}, a family of manifolds that generalizes the
classical Bott-Samelson cycles on $G/T$ [BS] (Lemma 3.4). In addition, 
\textsl{divided difference} in $K$-theory is introduced for \textsl{%
spherical represented involutions }(cf. \textbf{3.2}). In Section 4, by
resorting to the geometry of the adjoint representation, we interpret the
Bott-Samelson cycles on $G/T$ as certain twisted products of $2$-spheres,
and describe their $K$--rings in terms of the Cartan numbers of $G$ (Lemma
4.4). After determining the image of Demazure classes in the $K$--ring of a
Bott-Samelson cycle (Lemma 4.5), the theorems are established in Section 5.

\section{Preliminaries in topological $K$-theory}

All homologies (resp. cohomologies) will have integer coefficients unless
otherwise stated. If $f:M\rightarrow N$ is a continuous map between two
topological spaces $M$ and $N$, $f_{\ast }$ (resp. $f^{\ast }$) is the
homology (resp. cohomology) map induced by $f$, and $f^{!}:K(N)\rightarrow
K(M)$ is the induced map on the Grothendieck groups of topological complex
bundles. The involution on $K(M)$ by the complex conjugation is denoted by $%
\xi \rightarrow \overline{\xi }$, $\xi \in K(M)$.

Write $S^{2}$ for the $2$--dimensional sphere. If $M$ is an oriented closed
manifold $[M]\in H_{\dim M}(M)$ stands for the orientation class. The
Kronecker pairing, between cohomology and homology of a space $M$, will be
denoted by $<,>:H^{\ast }(M)\times H_{\ast }(M)\rightarrow \mathbb{Z}$.

Let $L_{\mathbb{C}}(M)$ be the set of isomorphism classes of complex line
bundles (i.e. of oriented real $2$--plane bundles) over $M$. The trivial
complex line bundle on $M$ is denoted by $1$. It is well known that

\begin{enumerate}
\item[(3.1)] sending a complex line bundle $\xi $ to its first Chern class $%
c_{1}(\xi )$ yields a one-to-one correspondence $c_{1}:L_{\mathbb{C}%
}(M)\rightarrow H^{2}(M)$ which is natural with respect to maps $%
M\rightarrow N$.
\end{enumerate}

\noindent Recall also from [AH] that

\begin{enumerate}
\item[(3.2)] if $M$ is a cell complex with even dimensional cells only, the 
\textsl{Chern character} $Ch:K(M)\rightarrow H^{\ast }(M;\mathbb{Q})$ is
injective.
\end{enumerate}

\textbf{3.1. }$S^{2}$\textbf{--bundle with a section. }Let $p:E\rightarrow M$
be a smooth, oriented $S^{2}$--bundle over an oriented manifold $M$ with a
section $s:M\rightarrow E$. As the normal bundle $\eta $ of the embedding $s$
is oriented by $p$ and has real dimension $2$, we may regard $\eta \in L_{%
\mathbb{C}}(M)$. We put $\xi =p^{!}(\eta )\in L_{\mathbb{C}}(E)$, $%
c=c_{1}(\xi )\in H^{2}(E)$.

The integral cohomology $H^{\ast }(E)$ can be described as follows. Denote
by $i:S^{2}\rightarrow E$ the fiber inclusion over a point $z\in M$, and
write $J:$ $E\rightarrow E$ for the involution given by the antipodal map in
each fiber sphere.

\begin{quote}
\textbf{Lemma 3.1 }(cf. [Du$_{2}$, Lemma 3.1])\textbf{. }\textsl{There
exists a unique class }$x\in H^{2}(E)$\textsl{\ such that }$s^{\ast
}(x)=0\in H^{\ast }(M)$\textsl{\ and }$<i^{\ast }(x),[S^{2}]>=1$\textsl{.
Furthermore }

\textsl{(1) }$H^{\ast }(E)$\textsl{, as a module over }$H^{\ast }(M)$\textsl{%
, has the basis }$\{1,x\}$\textsl{\ subject to the relation }$x^{2}+cx=0$%
\textsl{;}

\textsl{(2) the }$J^{\ast }$\textsl{\ acts} \textsl{identically on }$H^{\ast
}(M)\subset H^{\ast }(E)$ \textsl{and }$J^{\ast }(x)=-x-c$\textsl{.}
\end{quote}

Let $X\in L_{\mathbb{C}}(E)$ be with $c_{1}(X)=x\in H^{2}(E)$, where $x$ as
that in Lemma 3.1. Put $y=X-1\in K(E)$. The next result is seen as the $K$%
--theoretic analogous of Lemma 3.1, in which (1) is classical (cf. [A,
Proposition 2.5.3)).

\begin{quote}
\textbf{Lemma 3.2.} \textsl{Assume that }$M$\textsl{\ has a cell
decomposition with even dimensional cells only. Then}

\textsl{(1) }$K(E)$, \textsl{as a module over }$K(M)$\textsl{, has the basis 
}$\{1,y\}$\textsl{\ subject to the relation }$y^{2}=(\overline{\xi }-1)y$%
\textsl{. Furthermore}

\textsl{(2) }$\overline{y}=-\xi y$\textsl{;}

\textsl{(3) the }$J^{!}:K(E)\rightarrow K(E)$\textsl{\ acts identically\ on }%
$K(M)$\textsl{\ and satisfies }$J^{!}(y)=-y+(\overline{\xi }-1)$\textsl{.}
\end{quote}

\textbf{Proof.} Since $X$ restricts to the Hopf-line bundle on the fiber
sphere, $K(E)=K(M)[1,y]$. To verify the relation in (1) we compute

\begin{quote}
$Ch(y)=e^{x}-1=x+\frac{1}{2!}x^{2}+\frac{1}{3!}x^{3}+\cdots $

$\qquad \quad =-\frac{x}{c}(e^{-c}-1)$ (for $x^{n}=(-c)^{n-1}x$ by (1) of
Lemma 3.1).
\end{quote}

\noindent It follows that

\begin{quote}
$Ch(y^{2})=[\frac{x}{c}(e^{-c}-1)]^{2}=-\frac{x}{c}(e^{-c}-1)^{2}$ ($%
x^{2}=-cx$ by Lemma 3.1)

$\qquad \quad =Ch[(\overline{\xi }-1)y]$.
\end{quote}

\noindent This implies that $y^{2}=(\overline{\xi }-1)y$ by (3.2). This
completes (1).

In view of (1) we may assume that $\overline{y}(=\overline{X}-1)=a+by$, $%
a,b\in K(M)$. Multiplying both sides by $X=y+1$ yields

\begin{quote}
$\qquad -y=(a+by)(y+1)=a+(a+b\overline{\xi })y$ (by (1)).
\end{quote}

\noindent Coefficients comparison gives $a=0$; $a+b\overline{\xi }=-1$. This
shows (2).

Finally we show (3). From $J^{\ast }(x)=-x-c$ (by Lemma 3.1) and $%
c_{1}J^{!}=J^{\ast }c_{1}$ (by the naturality of (3.1)) we get $J^{!}(X)=%
\overline{X}\cdot \overline{\xi }$. From this one obtains

\begin{quote}
$J^{!}(y)=\overline{X}\cdot \overline{\xi }-1=(\overline{y}+1)\overline{\xi }%
-1=(-\xi y+1)\overline{\xi }-1$ (by (2))

$\qquad =-y+(\overline{\xi }-1)$. $\square $

\bigskip
\end{quote}

\textbf{3.2. Divided difference in }$K$\textbf{--theory. }A self-map $r$ of
a manifold $M$ is called \textsl{an involution} if $r^{2}=id:M\rightarrow M$%
. A $2$\textsl{--spherical representation} of the involution $(M;r)$ is a
system $f:(E;J)\rightarrow (M;r)$ in which $p:E\rightarrow M$ is an oriented 
$S^{2}$--bundle with a section $s$; and $f$ is a continuous map $%
E\rightarrow M$ that satisfies the following two constrains

\begin{enumerate}
\item[(3.3)] \noindent $\qquad f\circ s=id:M\rightarrow M$;$\qquad f\circ
J=r\circ f:$ $E\rightarrow M$,
\end{enumerate}

\noindent where $J$ is the involution on $E$ given by the antipodal map on
the fibers.

In view of the $K(M)$--module structure on $K(E)$ (cf. (1) of Lemma 3.2), a $%
2$--spherical representation (3.3) of the involution $(M,r)$ gives rise to
an additive operator $\Lambda _{f}:K(M)\rightarrow K(M)$ characterized
uniquely by (3.4) below, where the constraint $f\circ s=id$ implies that the
coefficient of $1$ is $p^{!}(z)$.

\textsl{The induced homomorphism }$f^{!}:K(M)\rightarrow K(E)$\textsl{\
satisfies}

\begin{enumerate}
\item[(3.4)] $\qquad \qquad \qquad f^{!}(z)=p^{!}(z)\cdot 1+p^{!}(\Lambda
_{f}(z))\cdot y$
\end{enumerate}

\noindent \textsl{for all }$z\in K(M)$\textsl{.}

The operator $\Lambda _{f}$ will be called \textsl{the divided difference}
associated to the $2$--spherical representation $f$ of the involution $(M;r)$%
.

\begin{quote}
\textbf{Lemma 3.3.} \textsl{Let }$\eta \in L_{\mathbb{C}}(M)$\textsl{\ the
normal bundle of the section }$s$\textsl{.} \textsl{Then}

\textsl{(1) }$r^{!}=Id+(\overline{\eta }-1)\Lambda _{f}$\textsl{\ }$:$%
\textsl{\ }$K(M)\rightarrow K(M)$\textsl{;}

\textsl{(2) }$\Lambda _{f}\circ r^{!}=-\Lambda _{f}$\textsl{.}
\end{quote}

\textbf{Proof.} Applying $J^{!}$ to (3.4) and using (3) of Lemma 3.2 to
write the resulting equality yield (note that $\xi =p^{!}(\eta )$ in Lemma
3.2)

\begin{center}
$\qquad J^{!}f^{!}(z)=p^{!}[z+(\overline{\eta }-1)\Lambda
_{f}(z)]-p^{!}[\Lambda _{f}(z)]y$.
\end{center}

\noindent On the other hand one gets from (3.4) that

\begin{center}
$\qquad f^{!}(r^{!}(z))=p^{!}(r^{!}(z))+p^{!}[\Lambda _{f}(r^{!}(z))]y$.
\end{center}

\noindent Since $J^{!}f^{!}=f^{!}r^{!}$ by (3.3), coefficients comparison
shows (1) and (2). $\square $

\bigskip

\textbf{3.3. The }$K$\textbf{-theory of a twisted product of }$S^{2}$. We
determine the $K$-rings for a class of manifolds specified below.

\begin{quote}
\textbf{Definition 5.} A smooth manifold $M$\ is called an oriented \textsl{%
twisted product of }$2$\textsl{--spheres\ of rank} $m$, written $M=\underset{%
1\leq i\leq m}{\propto }S^{2}$, if there is a tower of smooth maps
\end{quote}

\begin{center}
$M=M_{m}\overset{p_{m}}{\rightarrow }M_{m-1}\overset{p_{m-1}}{\rightarrow }%
\cdots \overset{p_{3}}{\rightarrow }M_{2}\overset{p_{2}}{\rightarrow }M_{1}%
\overset{p_{1}}{\rightarrow }M_{0}=\{z_{0}\}$
\end{center}

\begin{quote}
\noindent in which

1) $M_{0}$\ consists of a single point (as indicated);

2) each $p_{k}$\ is an oriented $S^{2}$--bundle with a fixed section $%
s_{k}:M_{k-1}\rightarrow M_{k}$.\ 
\end{quote}

Let $M=\underset{1\leq i\leq m}{\propto }S^{2}$ be a twisted product of $2$%
-spheres. Assign each $M_{k}$ with the base point $z_{k}=s_{k}\circ \cdots
\circ s_{1}(z_{0})\in M_{k}$ and denote by $h_{k}:$ $S^{2}\rightarrow M_{k}$
the inclusion of the fiber sphere of $p_{k}$ over the point $z_{k}$, $1\leq
k\leq m$. Consider the embedding $\iota _{k}:S^{2}\rightarrow M$ given by
the composition

\begin{center}
$s_{m}\circ \cdots \circ s_{k+1}\circ h_{k}:S^{2}\rightarrow
M_{k}\rightarrow M$.
\end{center}

\noindent Then the set $\{\iota _{1\ast }[S^{2}],\cdots ,\iota _{m\ast
}[S^{2}]\in H_{2}(M)\}$ of $2$-cycles forms a basis of $H_{2}(M)$ [Du$_{2}$,
Lemma 3.3]. Consequently, if we let $x_{i}\in H^{2}(M)$, $1\leq i\leq m$, be
the classes Kronecker dual to $\iota _{k\ast }[S^{2}]$ as $<x_{i},\iota
_{k\ast }[S^{2}]>=\delta _{ik}$, then

\begin{enumerate}
\item[(3.5)] the set $\{x_{1},\cdots ,x_{m}\}$ is a basis of $H^{2}(M)$ that
satisfies

$\qquad \qquad (s_{m}\circ \cdots \circ s_{k})^{\ast }(x_{k})=0$, $1\leq
k\leq m$.
\end{enumerate}

A set of numerical invariants for $M$ can now be extracted as follows. Let $%
\eta _{k}\in L_{\mathbb{C}}(M_{k-1})$ be the normal bundle of the embedding $%
s_{k}:M_{k-1}\rightarrow M_{k}$ with the induced orientation. Put $\xi
_{k}=(p_{k}\circ \cdots \circ p_{m})^{!}\eta _{k}\in L_{\mathbb{C}}(M)$. In
view of (3.5) we must have the expression in $H^{2}(M)$

$\qquad c_{1}(\xi _{1})=0$;

$\qquad c_{1}(\xi _{2})=a_{1,2}x_{1}$;

$\qquad c_{1}(\xi _{3})=a_{1,3}x_{1}+a_{2,3}x_{2};$

$\qquad \vdots $

$\qquad c_{1}(\xi _{m})=a_{1,m}x_{1}+\cdots +a_{m-1,m}x_{m-1}$

\noindent with $a_{i,j}\in \mathbb{Z}$.

\begin{quote}
\textbf{Definition 6.} With $a_{i,j}=0$\ for all $i\geq j$\ being
understood, the strictly upper triangular matrix $A=(a_{i,j})_{m\times m}$\
is called \textsl{the structure matrix of }$M=\underset{1\leq i\leq m}{%
\propto }S^{2}$\ relative to the basis $\{x_{1},\cdots ,x_{m}\}$ of $%
H^{2}(M) $.
\end{quote}

The ring $K(M)$ is determined by $A$ as follows. Let $X_{i}\in L_{\mathbb{C}%
}(M)$ be defined by $c_{1}(X_{i})=x_{i}$ (cf. (3.1), (3.5)) and set $%
y_{i}=X_{i}-1\in K(M)$. For a subset $I=[i_{1},\cdots ,i_{k}]\subseteq
\lbrack 1,\cdots ,m]$ we put

\begin{center}
$y_{I}=\{%
\begin{array}{c}
1\text{ if }k=0\text{ (i.e. }I=\{\emptyset \}\text{)}, \\ 
y_{I}=y_{i_{1}}\cdots y_{i_{k}}\text{ if }k\geq 1\text{.\ }%
\end{array}%
$
\end{center}

\noindent Let $q_{k},\overline{q_{k}}\in \mathbb{Z}[y_{1},\cdots ,y_{m}]$ be
defined in terms of $A$ as in Definition 3.

\begin{quote}
\textbf{Lemma 3.4. }\textsl{If }$M$\textsl{\ has the structure matrix }$%
A=(a_{i,j})_{m\times m}$\textsl{, then}

\textsl{(1) the set }$\{y_{I}\mid I\subseteq \lbrack 1,\cdots ,m]\}$\textsl{%
\ is a basis for }$K(M)$\textsl{;}

\textsl{(2) }$K(M)=\mathbb{Z}[y_{1},\cdots ,y_{m}]/<y_{k}^{2}=(\overline{%
q_{k}}-1)y_{k}$\textsl{, }$1\leq k\leq m>$\textsl{.}
\end{quote}

\textbf{Proof.} If $m=1$ then $M=S^{2}$ and it is clearly done. Assume next
that Lemma 3.4 holds for $m=n-1$, and consider now the case $m=n$.

Applying Lemma 3.2 to the $S^{2}$--bundle $M=M_{m}\overset{p_{m}}{%
\rightarrow }M_{m-1}$ we get

\begin{center}
$K(M)=K(M_{m-1})[1,y_{m}]$ ($y_{m}^{2}=(\overline{\xi _{m}}-1)y_{m}$).
\end{center}

\noindent This already shows (1).

For (2) it remains only to show $\overline{\xi _{m}}=\overline{q_{m}}\in 
\mathbb{Z}[y_{1},\cdots ,y_{m-1}]$. In view of (1) of Lemma 3.2 and by the
inductive hypothesis, we can assume that

\begin{enumerate}
\item[(3.6)] $\qquad \xi _{k}=q_{k}$,\quad\ $\overline{\xi _{k}}=\overline{%
q_{k}}$ for all $k\leq m-1$.
\end{enumerate}

\noindent From $c_{1}(\xi _{m})=a_{1,m}x_{1}+\cdots +a_{m-1,m}x_{m-1}$ we
find that

\begin{quote}
$\overline{\xi _{m}}=\prod_{a_{i,m}>0}\overline{X_{i}}^{a_{i,m}}%
\prod_{a_{i,m}<0}X_{i}^{-a_{i,m}}$

$\qquad =\prod_{a_{i,m}>0}(\overline{y_{i}}+1)^{a_{i,m}}%
\prod_{a_{i,m}<0}(y_{i}+1)^{-a_{i,m}}$ (since $y_{i}=X_{i}-1$)

$\qquad =\prod_{a_{i,m}>0}(-\xi
_{i}y_{i}+1)^{a_{i,m}}\prod_{a_{i,m}<0}(y_{i}+1)^{-a_{i,m}}$ (by Lemma 3.2)

$\qquad
=\prod_{a_{i,m}>0}(-q_{i}y_{i}+1)^{a_{i,m}}%
\prod_{a_{i,m}<0}(y_{i}+1)^{-a_{i,m}}$ (by (3.6))

$\qquad =\overline{q_{m}}$ (cf. Definition 3).
\end{quote}

\noindent Similarly, we have $\xi
_{m}=\prod_{a_{i,m}>0}X^{a_{i,m}}\prod_{a_{i,m}<0}\overline{X_{i}}%
^{-a_{i,m}}=q_{m}$.$\square $

\bigskip

According to (1) of Lemma 3.4, any polynomial $f\in \mathbb{Z}[y_{1},\cdots
,y_{m}]$ can be expanded (uniquely) as a linear combination of the $y_{I}$

\begin{center}
$f=\sum a_{I}(f)y_{I}$, $I\subseteq \lbrack 1,\cdots ,m]$,
\end{center}

\noindent where the correspondences $a_{I}:\mathbb{Z}[y_{1},\cdots
,y_{m}]\rightarrow \mathbb{Z}$ by $f\rightarrow a_{I}(f)$ are clearly
additive. Indeed, problems concerning computing in the $K(M)$ ask an
effective algorithm to evaluate the $a_{I}$. The case $I=[1,\cdots ,m]$ will
be relevant to us and whose solution brings us the operator $\Delta _{A}$
given by Definition 4.

\begin{quote}
\textbf{Lemma 3.5. }\textsl{If }$M$\textsl{\ has the structure matrix }$%
A=(a_{i,j})_{m\times m}$\textsl{, then}
\end{quote}

\begin{center}
$a_{[1,\cdots ,m]}=\Delta _{A}\circ \quad _{(m)}:\mathbb{Z}[y_{1},\cdots
,y_{m}]\overset{(m)}{\rightarrow }\mathbb{Z}[y_{1},\cdots
,y_{m}]_{(m)}\rightarrow \mathbb{Z}$\textsl{. }
\end{center}

\begin{quote}
\noindent \textsl{In particular, }$a_{[1,\cdots ,m]}(f)=0$\textsl{\ if }$%
f_{(m)}=0$\textsl{.}
\end{quote}

\textbf{Proof.} This is parallel to the proof of Proposition 2 in [Du$_{1}$].%
$\square $

\section{Bott-Samelson cycles and Demazure classes}

With respect to the fixed regular point $\alpha \in L(T)$ the adjoint
representation $Ad:G\rightarrow L(G)$ gives rise to a smooth embedding

\begin{center}
$\varphi :G/T\rightarrow L(G)\quad \quad $by $\varphi (gT)=Ad_{g}(\alpha )$.
\end{center}

\noindent In this way $G/T$ becomes a submanifold of the Euclidean space $%
L(G)$. By resorting to the geometry of this embedding we recover the
Demazure operators on $K(G/T)$ in \textbf{4.3}, and\textbf{\ }the classical
Bott-Samelson cycle $\varphi _{0,\beta _{1},\cdots ,\beta _{k}}:S(\alpha
;\beta _{1},\cdots ,\beta _{k})\rightarrow G/T$ (associated to a sequence of
simple roots $\beta _{1},\cdots ,\beta _{k}$) on $G/T$ in \textbf{4.4}. As
application of Lemma 3.4, the ring $K(S(\alpha ;\beta _{1},\cdots ,\beta
_{k}))$ is described in terms of the Cartan numbers of $G$ (cf. Lemma 4.4).
The main result in this section is Lemma 4.5, which specifies the $\varphi
_{0,\beta _{1},\cdots ,\beta _{k}}^{!}$--image of a Demazure class in $%
K(S(\alpha ;\beta _{1},\cdots ,\beta _{k}))$.\ \ 

\textbf{4.1.} \textbf{Geometry from the adjoint representation.} Let $\Phi
^{+}\subset L(T)$ (resp. $\Delta \subset L(T)$) be the set of positive
(resp. simple) roots relative to $\alpha $ ([Hu, p.35]). Assume that the
Cartan decomposition of the Lie algebra $L(G)$ relative to $T\subset G$ is

\begin{center}
$L(G)=L(T)\oplus _{\beta \in \Phi ^{+}}F_{\beta }$,
\end{center}

\noindent where $F_{\beta }$ is the root space, viewed as a real $2$-plane,
belonging to the root $\beta \in \Phi ^{+}$ ([Hu, p.35]). It is known (cf.
[HPT,p.426-427; or Du$_{2}$, Sect.4]) that

\begin{enumerate}
\item[(4.1)] The subspaces $\oplus _{\beta \in \Phi ^{+}}F_{\beta }$ and $%
L(T)$ of $L(G)$ are tangent and normal to $G/T$ at $\alpha $ respectively;

\item[(4.2)] The tangent bundle to $G/T$ has a canonical orthogonal
decomposition into the sum of integrable $2$-plane bundles $\oplus _{\beta
\in \Phi ^{+}}E_{\beta }$ with $E_{\beta }(\alpha )=F_{\beta }$.

\item[(4.3)] The leaf of the integrable subbundle $E_{\beta }$ through a
point $z\in G/T$, denoted by $S(z;\beta )$, is a $2$-sphere that carries a
preferred orientation.

\item[(4.4)] Via the embedding $\varphi $, the canonical action of $W$ on $%
G/T$ can be given in terms of the $W$ action on $L(T)$ as $%
w(z)=Ad_{g}(w(\alpha ))$ if $z=Ad_{g}(\alpha )\in G/T$, $w\in W$ (cf. [BS]).
\end{enumerate}

\textbf{4.2. Demazure basis of }$K(G/T)$. For a complete subvariety $%
Y\subset G/T$ the \textsl{Euler--Poincar\'{e} characteristic} relative to $Y$
is the homomorphism $\chi (Y,-):K(G/T)\rightarrow \mathbb{Z}$ defined by

\begin{center}
$[\digamma ]\rightarrow \chi (Y,\digamma
)=\sum\limits_{j}(-1)^{j}h^{j}(\digamma \mid Y)$,
\end{center}

\noindent where $\digamma \mid Y$ means the restriction of $\digamma $ on $Y$%
, and where $h^{j}(\digamma \mid Y)$ denotes the dimension of the $j^{th}$
cohomology group of $h^{j}(\digamma \mid Y)$. The following characterization
of Demazure basis is due to B. Kostant and S. Kumar (compare [KK, (3.39)
Proposition] with [D, Proposition 7]).

\begin{quote}
\textbf{Definition 7. }The \textsl{Demazure basis} $\{a_{w}\in K(G/T)\mid
w\in W\}$ of the ring $K(G/T)$ is defined by
\end{quote}

\begin{center}
$\chi (X_{w},a_{u})=\delta _{w,u}$, $w,u\in W$,
\end{center}

\begin{quote}
where $X_{w}$ is the Schubert class on $G/T$ associated to $w$.
\end{quote}

\bigskip

\textbf{4.3. Divided difference on }$K(G/T)$\textbf{\ associated to a root. }%
Each root\textbf{\ }$\beta \in \Phi ^{+}$ gives rise to an involution $%
r_{\beta }:G/T\rightarrow G/T$ in the fashion of (4.4), and defines also the
subspace

\begin{center}
$S(\beta )=\{(z,z_{1})\in G/T\times G/T\mid z_{1}\in S(z;\beta )\}$
\end{center}

\noindent in view of (4.3). Projection $p_{\beta }:S(\beta )\rightarrow G/T$
onto the first factor is easily seen to be a $S^{2}$--bundle (with $%
S(z;\beta )$ as the fiber over $z\in G/T$). The map $s_{\beta
}:G/T\rightarrow S(\beta )$ by $s_{\beta }(z)=(z,z)$ furnishes $p_{\beta }$
with a ready-made section.

Let $x\in H^{2}(S(\beta ))$ be such that $s_{\beta }^{\ast }(x)=0$ and $%
<i^{\ast }(x),[S(z;\beta )]>=1$\ (cf. Lemma 3.1) and set $y=X-1\in K(S(\beta
))$, where $X\in L_{\mathbb{C}}(S(\beta ))$ is defined by $c_{1}(X)=x$ (cf.
(3.1)). Since the normal bundle of the embedding $s_{\beta }$ is easily seen
to be $E_{\beta }\in L_{\mathbb{C}}(G/T)$, one infers from Lemma 3.2 that

\begin{quote}
\textbf{Lemma 4.1.}\textsl{\ }$K(S(\beta ))=K(G/T)[1,y]/<y^{2}=(\overline{%
E_{\beta }}-1)y>$.
\end{quote}

Let $J_{\beta }$ be the involution on $S(\beta )$ given by the antipodal
maps in the fiber spheres, and let $f_{\beta }:S(\beta )\rightarrow G/T$ be
the projection onto the second factor. Then, as is clear,

\begin{center}
$f_{\beta }\circ s_{\beta }=id:G/T\rightarrow G/T$;$\quad f_{\beta }\circ
J_{\beta }=r_{\beta }\circ f_{\beta }:S(\beta )\rightarrow G/T$.
\end{center}

\noindent That is, \textsl{the map }$f_{\beta }:(S(\beta ),J_{\beta
})\rightarrow (G/T,r_{\beta })$\textsl{\ is a }$2$\textsl{-spherical
representation of the involution }$(G/T,r_{\beta })$\textsl{\ }(cf. \textbf{%
3.2}).

Abbreviate the divided difference $\Lambda _{f_{\beta }}:K(G/T)\rightarrow
K(G/T)$ associated to $f_{\beta }$ by $\Lambda _{\beta }$. The next result,
essentially due to Kostant and Kumar [KK], is the key in the proof of Lemma
4.5.

\begin{quote}
\textbf{Lemma 4.2.} Let $\{a_{w}\in K(G/T)\mid w\in W\}$ be the \textsl{%
Demazure basis} of $K(G/T)$, and let $\beta \in \Delta $ be a simple root.
Then
\end{quote}

\begin{center}
$\Lambda _{\beta }(\overline{a_{w}})=\{%
\begin{array}{c}
E_{\beta }\cdot \overline{a_{w\cdot r_{\beta }}}\text{\quad if }%
l(w)>l(wr_{\beta })\text{;} \\ 
-E_{\beta }\cdot \overline{a_{w}}\text{\quad if }l(w)<l(wr_{\beta })\text{.}%
\end{array}%
$
\end{center}

\textbf{Proof.} Recall from [KK,PR$_{1}$,D] that the classical \textsl{%
Demazure operator} $T_{\beta }:K(G/T)\rightarrow K(G/T)$ is given by

\begin{center}
$T_{\beta }(u)=$ $\frac{E_{\beta }\cdot u-r_{\beta }^{!}(u)}{E_{\beta }-1}$.
\end{center}

\noindent Substituting in $r_{\beta }^{!}=Id+(\overline{E_{\beta }}%
-1)\Lambda _{\beta }$ (lemma 3.3) yields

\begin{center}
$T_{\beta }(u)=u+\overline{E_{\beta }}\Lambda _{\beta }(u)$.
\end{center}

\noindent That is $\Lambda _{\beta }=E_{\beta }(T_{\beta
}-Id):K(G/T)\rightarrow K(G/T)$.

On the other hand combining [KK, Proposition (2.22),(d)] with [KK,
Proposition (3.39)] one gets

\begin{center}
$T_{\beta }(\overline{a_{w}})=\{%
\begin{array}{c}
\overline{a_{w}}+\overline{a_{wr_{\beta }}}\text{\quad if }l(w)>l(wr_{\beta
})\text{;} \\ 
0\text{, otherwise.\ \qquad \quad \qquad \qquad }%
\end{array}%
$
\end{center}

\noindent This completes the proof.$\square $

\bigskip

\textbf{4.4.} \textbf{Bott-Samelson cycles and their }$K$\textbf{-rings. }%
Given an ordered sequence $(\beta _{1},\cdots ,\beta _{k})$ of simple roots
(repetitions like $\beta _{i}=\beta _{j}$ for some $1\leq i<j\leq k$ may
occur), we set

\begin{center}
$S(\alpha ;\beta _{1},\cdots ,\beta _{k})=\{(z_{0},z_{1},\cdots ,z_{k})\in
G/T\times \cdots \times G/T\mid z_{0}=\alpha $; $z_{i}\in S(z_{i-1};\beta
_{i})\}$.
\end{center}

\noindent It is furnished with the structure of oriented twisted product of $%
2$-spheres of rank $k$ by the maps

\begin{center}
$S(\alpha ;\beta _{1},\cdots ,\beta _{i})\overset{p_{i}}{\underset{s_{i}}{%
\rightleftarrows }}S(\alpha ;\beta _{1},\cdots ,\beta _{i-1})$,

$p_{i}(z_{0},\cdots ,z_{i})=(z_{0},\cdots ,z_{i-1})$; $s_{i}(z_{0},\cdots
,z_{i-1})=(z_{0},\cdots ,z_{i-1},z_{i-1})$.
\end{center}

\noindent One has also the ready-made maps

\begin{center}
$\varphi _{0,\beta _{1},\cdots ,\beta _{k}}:S(\alpha ;\beta _{1},\cdots
,\beta _{k})\rightarrow G/T$ by $(z_{0},\cdots ,z_{k})\rightarrow z_{k}$

$\widehat{\varphi }_{0,\beta _{1},\cdots ,\beta _{k}}:S(\alpha ;\beta
_{1},\cdots ,\beta _{k})\rightarrow S(\beta _{k})$ by $(z_{0},\cdots
,z_{k})\rightarrow (z_{k-1},z_{k})$
\end{center}

\noindent that clearly satisfy

\begin{enumerate}
\item[(4.5)] \noindent\ $\varphi _{0,\beta _{1},\cdots ,\beta _{k}}=f_{\beta
_{k}}\circ \widehat{\varphi }_{0,\beta _{1},\cdots ,\beta _{k-1}}:S(\alpha
;\beta _{1},\cdots ,\beta _{k})\rightarrow G/T$;

\item[(4.6)] the commutative diagrams
\end{enumerate}

\begin{center}
$%
\begin{array}{ccccc}
S(\alpha ;\beta _{1},\cdots ,\beta _{k}) & \overset{\widehat{\varphi }%
_{0,\beta _{1},\cdots ,\beta _{k-1}}}{\rightarrow } & S(\beta _{k}) &  &  \\ 
p_{k-1}\downarrow \uparrow s_{k-1} &  & p_{\beta _{k}}\downarrow \uparrow
s_{\beta _{k}} & \overset{f_{\beta _{k}}}{\searrow } &  \\ 
S(\alpha ;\beta _{1},\cdots ,\beta _{k-1}) & \overset{\varphi _{0,\beta
_{1},\cdots ,\beta _{k-1}}}{\rightarrow } & G/T & \overset{f_{\beta
_{k}}\circ s_{\beta _{k}}=id}{\rightarrow } & G/T%
\end{array}%
$
\end{center}

\noindent in which $\widehat{\varphi }_{0,\beta _{1},\cdots ,\beta _{k-1}}$
is a bundle map over $\varphi _{0,\beta _{1},\cdots ,\beta _{k-1}}$.

\begin{quote}
\textbf{Definition 8 }([Du$_{2}$, \textbf{7.2]}). \textsl{The map (4.5) is
called the Bott-Samelson cycle associated to the sequence }$\beta
_{1},\cdots ,\beta _{k}$\textsl{\ of simple roots.}
\end{quote}

Let $\iota _{i}:S(\alpha ,\beta _{i})\rightarrow S(\alpha ;\beta _{1},\cdots
,\beta _{k})$ be the embedding specified by

\begin{center}
$\iota _{i}(\alpha ,z^{\prime })=(z_{0},\cdots ,z_{k})$,
\end{center}

\noindent where $z_{0}=\cdots =z_{i-1}=\alpha $, $z_{i}=\cdots
=z_{k}=z^{\prime }$. Then the cycles $\iota _{i\ast }[S(\alpha ,\beta
_{i})]\in H_{2}(S(\alpha ;\beta _{1},\cdots ,\beta _{k}))$, $1\leq i\leq k$,
form a basis of $H_{2}(S(\alpha ;\beta _{1},\cdots ,\beta _{k}))$ (cf. 
\textbf{3.3}). Let $x_{i}\in H^{2}(S(\alpha ;\beta _{1},\cdots ,\beta _{k}))$
be the basis Kronecker dual to the $\iota _{j\ast }[S(\alpha ,\beta _{j})]$
as $<x_{i},\iota _{j\ast }[S(\alpha ,\beta _{j})]>=\delta _{ij}$. The next
result was shown in [Du$_{2}$, Lemma 4.5].

\begin{quote}
\textbf{Lemma 4.3.} \textsl{The structure matrix }$A=(a_{i,j})_{k\times k}$ 
\textsl{of }$S(\alpha ;\beta _{1},\cdots ,\beta _{k})$ \textsl{relative to} $%
\{x_{1},\cdots ,x_{k}\}$ \textsl{is given by the Cartan numbers of }$G$%
\textsl{\ as}
\end{quote}

\begin{center}
$a_{i,j}=\{%
\begin{array}{c}
\beta _{j}\circ \beta _{i}\ \quad \text{\textsl{if} }i<j\text{;} \\ 
0\ \quad \text{\textsl{if} }i\geq j\text{\textsl{.}\quad \quad }%
\end{array}%
$
\end{center}

Let $X_{i}\in L_{\mathbb{C}}(S(\alpha ;\beta _{1},\cdots ,\beta _{k}))$ be
defined by $c_{1}(X_{i})=x_{i}$. Set $y_{i}=X_{i}-1\in K(S(\alpha ;\beta
_{1},\cdots ,\beta _{k}))$. Let $\overline{q_{k}}\in \mathbb{Z}[y_{1},\cdots
,y_{m}]$ be defined in terms of $A$ as in Definition 3. Combining Lemma 3.4,
Lemma 3.5 with Lemma 4.3 yields the next result.

\begin{quote}
\textbf{Lemma 4.4. }\textsl{Let }$M=S(\alpha ;\beta _{1},\cdots ,\beta _{k})$%
\textsl{. Then}

\textsl{(1) the set }$\{y_{I}\mid I\subseteq \lbrack 1,\cdots ,m]\}$\textsl{%
\ is a basis of }$K(M)$\textsl{;}

\textsl{(2) }$K(M)=\mathbb{Z}[y_{1},\cdots ,y_{k}]/<y_{i}^{2}=(\overline{%
q_{i}}-1)y_{i}$\textsl{, }$1\leq i\leq k>$\textsl{;}

\textsl{(3) if }$f\in \mathbb{Z}[y_{1},\cdots ,y_{k}]$\textsl{\ with }$%
f=\sum a_{I}(f)y_{I}$\textsl{, then }
\end{quote}

\begin{center}
$a_{[1,\cdots ,k]}(f)=\Delta _{A}(f_{(k)})$\textsl{. }
\end{center}

\begin{quote}
\textsl{In particular, }$a_{[1,\cdots ,k]}(f)=0$\textsl{\ if} \textsl{\ }$%
f_{(k)}=0$\textsl{.}
\end{quote}

\textbf{4.5. The induced map of a Bott-Samelson cycle. }Given a sequence $%
\beta _{1},\cdots ,\beta _{k}$ of simple roots consider the induced ring map

\begin{center}
$\varphi _{0,\beta _{1},\cdots ,\beta _{k}}^{!}:K(G/T)\rightarrow K(S(\alpha
;\beta _{1},\cdots ,\beta _{k}))$.
\end{center}

\noindent The ring $K(S(\alpha ;\beta _{1},\cdots ,\beta _{k}))$ has the
additive basis $\{y_{I}\mid $\textsl{\ }$I\subseteq $\textsl{\ }$[1,\cdots
,k]\}$ by Lemma 4.4. \textsl{\ }Let $\{a_{w}\mid w\in W\}$\ be the Demazure
basis of $K(G/T)$.

\begin{quote}
\textbf{Lemma 4.5.}\textsl{\ The induced map} $\varphi _{0,\beta _{1},\cdots
,\beta _{k}}^{!}$ \textsl{is given by}
\end{quote}

\begin{center}
$[\varphi _{0,\beta _{1},\cdots ,\beta
_{k}}]^{!}(a_{w})=(-1)^{l(w)}\sum\limits_{I\subseteq \lbrack 1,\cdots
,k],\beta (I)\thicksim w}y_{I}$\textsl{.}
\end{center}

\textbf{Proof.} It suffices to show that

\begin{enumerate}
\item[(4.7)] $\qquad \lbrack \varphi _{0,\beta _{1},\cdots ,\beta _{k}}]^{!}(%
\overline{a_{w}})=(-1)^{l(w)}\sum\limits_{I\subseteq \lbrack 1,\cdots
,k],\beta (I)\thicksim w}\overline{y}_{I}$\textsl{,}
\end{enumerate}

\noindent for, the complex conjugation of (4.7) yields the Lemma. To this
end we compute

\begin{quote}
$\varphi _{0,\beta _{1},\cdots ,\beta _{k}}^{!}(\overline{a_{w}})=\widehat{%
\varphi }_{0,\beta _{1},\cdots ,\beta _{k-1}}^{!}(f_{\beta _{k}}^{!}((%
\overline{a_{w}})))$ (by (4.5))

\noindent $\quad =\widehat{\varphi }_{o,\beta _{1},\cdots ,\beta
_{k-1}}^{!}[p_{\beta _{k}}^{!}(\overline{a_{w}})+p_{\beta _{k}}^{!}(\Lambda
_{\beta _{k}}(\overline{a_{w}}))y_{k}]$ (by (3.4))

\noindent $\quad =\varphi _{0,\beta _{1},\cdots ,\beta _{k-1}}^{!}(\overline{%
a_{w}})+\varphi _{0,\beta _{1},\cdots ,\beta _{k-1}}^{!}(\Lambda _{\beta
_{k}}(\overline{a_{w}}))y_{k}$ (by (4.6))

\noindent $\quad =\{%
\begin{array}{c}
\varphi _{0,\beta _{1},\cdots ,\beta _{k-1}}^{!}(\overline{a_{w}})+\varphi
_{0,\beta _{1},\cdots ,\beta _{k-1}}^{!}(\overline{a_{wr_{\beta _{k}}}}%
)E_{\beta _{k}}\cdot y_{k}\text{\quad if }l(w)>l(wr_{\beta })\text{;} \\ 
\varphi _{0.\beta _{1},\cdots ,\beta _{k-1}}^{!}(\overline{a_{w}})-\varphi
_{0,\beta _{1},\cdots ,\beta _{k-1}}^{!}(\overline{a_{w}})E_{\beta
_{k}}\cdot y_{k}\text{,\quad otherwise,\qquad \quad \quad }%
\end{array}%
$
\end{quote}

\noindent where the last equality is by Lemma 4.2. From $\overline{y_{k}}%
=-E_{\beta _{k}}\cdot y_{k}$ (Lemma 3.2) we obtain

\begin{center}
\noindent $\varphi _{0,\beta _{1},\cdots ,\beta _{k}}^{!}(\overline{a_{w}}%
)=\{%
\begin{array}{c}
\varphi _{0,\beta _{1},\cdots ,\beta _{k-1}}^{!}(\overline{a_{w}})-\varphi
_{0,\beta _{1},\cdots ,\beta _{k-1}}^{!}(\overline{a_{wr_{\beta _{k}}}})%
\overline{y_{k}}\text{ if }l(w)>l(wr_{\beta })\text{;} \\ 
\varphi _{0,\beta _{1},\cdots ,\beta _{k-1}}^{!}(\overline{a_{w}})+\varphi
_{0,\beta _{1},\cdots ,\beta _{k-1}}^{!}(\overline{a_{w}})\overline{y_{k}}%
\text{, otherwise.\quad \qquad \quad }%
\end{array}%
$
\end{center}

\noindent This, together with an easy induction on $k$, reduces the proof of
(4.7) to the general properties of the Demazure classes. Let $f:X\rightarrow
G/T$ be a continuous map from a connected $2m$--dimensional $CW$-complex $X$%
. Then

(1) $f^{!}(\overline{a_{u}})=0$ whenever $l(u)>m$;

(2) $f^{!}(\overline{a_{e}})=1$ if $X$ is a single point, where $e\in W$ is
the group unit,

\noindent Indeed, (1) follows from Assertion III in [KK, (3.26)] and [KK,
(3.39)]; (2) can be deduced from (1) and

(i) $\Omega _{e}=1$ (cf. [PR$_{2}$, Corollary 2.5] and the footnote in 
\textbf{5.2});

(ii) $\Omega _{e}=\sum_{w\in W}a_{w}$ (cf. Lemma 5.1 below).$\square $

\begin{quote}
\textbf{Corollary 1.} \textsl{Let }$e\in W$\textsl{\ be the group unit. Then}

\textsl{(1)} $[\varphi _{0,\beta _{1},\cdots ,\beta
_{k}}]^{!}(a_{e})=\prod\limits_{1\leq i\leq k}$ $(1+y_{i})$\ 

\noindent \textsl{and for} $w\neq e$

\textsl{(2)} $[\varphi _{0,\beta _{1},\cdots ,\beta _{k}}]^{!}(a_{w})=\{%
\begin{array}{c}
0\text{\quad \textsl{if} }l(w)>k\text{\textsl{;}\qquad \quad \qquad \qquad
\qquad } \\ 
(-1)^{k}\delta _{w,r_{\beta _{1}}\cdots r_{\beta _{k}}}y_{[1,\cdots ,m]}%
\text{\quad \textsl{if} }l(w)=k\text{\textsl{.}}%
\end{array}%
$
\end{quote}

\textbf{Example 1. }Let $\Delta =\{\alpha _{1},\alpha _{2}\}$ be a set of
simple roots of $G_{2}$ in which $\alpha _{1}$ is the short root [Hu, p.57].
The Weyl group $W$ of $G_{2}$ is generated by $\sigma _{i}$, $i=1,2$, the
reflection in the hyperplane $L_{i}\subset L(T)$ perpendicular to $\alpha
_{i}$. If we take $u=\sigma _{1}\sigma _{2}\sigma _{1}\sigma _{2}\sigma _{1}$%
, then Lemma 4.5 yields that

\begin{quote}
(a) $[\varphi _{0,(\alpha _{1},\alpha _{2},\alpha _{1},\alpha _{2},\alpha
_{1})}]^{!}(a_{u})=x_{1}x_{2}x_{3}x_{4}x_{5}$;

(b) $[\varphi _{0,(\alpha _{2},\alpha _{1},\alpha _{2},\alpha _{1},\alpha
_{2})}]^{!}(a_{u})=0$;

(c) $[\varphi _{0,(\alpha _{1},\alpha _{2},\alpha _{1},\alpha _{2},\alpha
_{1},\alpha _{2})}]^{!}(a_{u})=x_{1}x_{2}x_{3}x_{4}x_{5}(1+x_{6})$.
\end{quote}

\section{Multiplication in the ring $K(G/H)$}

Lemma 4.5 enables us to establish Theorem 1 (resp. Theorem 2) by computation
in the simpler ring $K(S(\alpha ;\beta _{1},\cdots ,\beta _{k}))$ (cf. Lemma
4.4).

\bigskip

\textbf{5.1. Proof of Theorem 1} (cf. \S 2)\textbf{. }Let $w=r_{\beta
_{1}}\cdot \cdots \cdot r_{\beta _{m}}$, $\beta _{i}\in \Delta $ be a
reduced decomposition of a $w\in W$, and let $A_{w}=(a_{i,j})_{m\times m}$
be the associated Cartan matrix (cf. Definition 1). Consider the
Bott-Samelson cycle $\varphi _{0,\beta _{1},\cdots ,\beta _{m}}:S(\alpha
;\beta _{1},\cdots ,\beta _{m})\rightarrow G/T$ associated to the sequence $%
\beta _{1},\cdots ,\beta _{m}$ of simple roots. Applying the induced ring
map $\varphi _{0,\beta _{1},\cdots ,\beta _{m}}^{!}$ to (1.1) yields in $%
K(S(\alpha ;\beta _{1},\cdots ,\beta _{m}))$ that

\begin{quote}
$\qquad \qquad \varphi _{0,\beta _{1},\cdots ,\beta _{m}}^{!}[a_{u}\cdot
a_{v}]=\sum\limits_{x\in W}C_{u,v}^{x}\varphi _{0,\beta _{1},\cdots ,\beta
_{m}}^{!}[a_{x}]$

$\qquad =(-1)^{l(w)}C_{u,v}^{w}y_{[1,\cdots ,m]}+\sum_{l(x)\leq
l(w)-1}C_{u,v}^{x}\varphi _{0,\beta _{1},\cdots ,\beta _{m}}^{!}[a_{x}]$,
\end{quote}

\noindent where the second equality follows from (2) of Corollary 1. Using
Lemma 4.5 to rewrite this equation yields

\begin{quote}
\bigskip $\qquad (-1)^{l(u)+l(v)}[\sum\limits_{\beta (L)\thicksim
u}y_{L}][\sum\limits_{\beta (K)\thicksim v}y_{K})]$

$=(-1)^{l(w)}C_{u,v}^{w}y_{[1,\cdots ,m]}+\sum_{l(x)\leq
l(w)-1}(-1)^{l(x)}C_{u,v}^{x}[\sum\limits_{\beta (J)\thicksim x}y_{J}]$,
\end{quote}

\noindent where $L,K,J\subseteq \lbrack 1,\cdots ,m]$. Finally, comparing
the coefficients of the monomial $y_{[1,\cdots ,m]}$ on both sides by using
(3) of Lemma 4.4, we obtain

\begin{quote}
\bigskip $\qquad \qquad (-1)^{l(u)+l(v)}\Delta _{A_{w}}[(\sum\limits_{\beta
(L)\thicksim u}y_{L})(\sum\limits_{\beta (K)\thicksim v}y_{K})]_{(m)}$

$\qquad =(-1)^{l(w)}C_{u,v}^{w}+\sum\limits_{\substack{ l(u)+l(v)\leq
l(x)\leq l(w)-1  \\ (\beta _{1},\cdots ,\beta _{m})\thicksim x}}%
(-1)^{l(x)}C_{u,v}^{x}$.
\end{quote}

\noindent This finishes the proof.$\square $

\bigskip

\textbf{Example 2. }Continuing from Example 1 we take

\begin{center}
$u=\sigma _{1}\sigma _{2}\sigma _{1}\sigma _{2}\sigma _{1}$;$\quad v=\sigma
_{2}\sigma _{1}\sigma _{2}\sigma _{1}\sigma _{2}$;$\quad w=\sigma _{1}\sigma
_{2}\sigma _{1}\sigma _{2}\sigma _{1}\sigma _{2}$.
\end{center}

\noindent Then the $u,v\in W$ are the only elements of length $5$, and $w$
is the element of highest length. Applying Theorem 1 we compute the
structure constants appearing in the expansion

\begin{center}
$a_{e}a_{u}=C_{e,u}^{u}a_{u}+C_{e,u}^{v}a_{v}+C_{e,u}^{w}a_{w}$.
\end{center}

\noindent In views of (a) in Example 1 and (1) in Corollary 1 we have

\begin{quote}
$C_{e,u}^{u}=\Delta _{A_{u}}[x_{1}x_{2}x_{3}x_{4}x_{5}\cdot \prod_{1\leq
i\leq 5}(1+x_{i})]_{(5)}$

$\qquad =\Delta _{A_{u}}(x_{1}x_{2}x_{3}x_{4}x_{5})=1$.
\end{quote}

\noindent Similarly, we get

\begin{quote}
$C_{e,u}^{v}=\Delta _{A_{v}}[0\cdot \prod_{1\leq i\leq 5}(1+x_{i})]_{(5)}=0$;

$C_{e,u}^{w}=(-1)^{5}\Delta _{A_{w}}[x_{1}x_{2}x_{3}x_{4}x_{5}(1+x_{6})\cdot
\prod_{1\leq i\leq 5}(1+x_{i})]_{(6)}$

$\qquad \qquad -[(-1)^{5}(C_{e,u}^{u}+C_{e,u}^{v})]$

$\qquad =-\Delta _{A_{w}}[x_{1}x_{2}x_{3}x_{4}x_{5}(x_{1}+\cdots
+x_{5}+2x_{6})]+1$

$\qquad =-2+1=-1$,
\end{quote}

\noindent where the three $\Delta _{A}(f)$'s concerned in the above
computation are directly evaluated from Definition 4 without resorting to
the specialities of $A$, thanks to the simplicities of the polynomials $f$
involved. Summarizing we get, in the ring $K(G_{2}/T)$, that

\begin{center}
$a_{e}a_{u}=a_{u}-a_{w}$.
\end{center}

\bigskip

\textbf{5.2.} The method establishing Theorem 1 is directly applicable to
find a formula for the structure constants $K_{u,v}^{w}(H)$ for multiplying
Grothendieck classes in the $K(G/H)$.

We begin with the simpler case $H=T$ (a maximal torus in $G$). Abbreviate $%
X_{w}(H)$ by $X_{w}$, $\Omega _{w}(H)$ by $\Omega _{w}$\footnote{%
The $\Omega _{w}$ corresponds to $\mathcal{O}_{w_{0}w}$ \ in [Br$_{2}$,KK,PR$%
_{2}$].} and $K_{u,v}^{w}(H)$ by $K_{u,v}^{w}$. The transition between the
two bases $\{a_{w}\mid w\in W\}$ and $\{\Omega _{w}\mid w\in W\}$ of $K(G/T)$
has been determined by Kostant and Kumar in [KK, Proposition 4.13; 3.39].

\begin{quote}
\textbf{Lemma 5.1.} \textsl{In the ring }$K(G/T)$\textsl{\ one has }
\end{quote}

\begin{center}
$\Omega _{w}=\sum\limits_{w\leq u}a_{u}$\textsl{,\qquad }$%
a_{w}=\sum\limits_{w\leq u}(-1)^{l(u)-l(w)}\Omega _{u}$
\end{center}

\begin{quote}
\noindent \textsl{where }$w\leq u$\textsl{\ means }$X_{w}\subseteq X_{u}$%
\textsl{.}
\end{quote}

Combining Lemma 5.1 with Lemma 4.5 gives the next result.

\begin{quote}
\textbf{Lemma 5.2.} \textsl{Let }$\beta _{1},\cdots ,\beta _{k}$ \textsl{be
a sequence of simple roots. With respect to the Grothendieck basis the
induced map} $\varphi _{0,\beta _{1},\cdots ,\beta _{k}}^{!}$\textsl{\ is
given by}
\end{quote}

\begin{center}
$[\varphi _{0,\beta _{1},\cdots ,\beta _{k}}]^{!}(\Omega
_{w})=\sum\limits_{I\subseteq \lbrack 1,\cdots ,k]}b_{I}(w)y_{I}$\textsl{,}
\end{center}

\begin{quote}
\noindent \textsl{where }$b_{I}(w)=\sum\limits_{\beta (I)\sim u,u\geq
w}(-1)^{l(u)}$\textsl{.}
\end{quote}

\bigskip

Based on Lemma 5.2, an argument parallel to the proof of Theorem 1 yields

\begin{quote}
\textbf{Theorem 2.} \textsl{Assume that }$w=r_{\beta _{1}}\cdot \cdots \cdot
r_{\beta _{m}}$\textsl{, }$\beta _{i}\in \Delta $\textsl{, is a reduced
decomposition of a }$w\in W$\textsl{, and let }$A_{w}=(a_{i,j})_{m\times m}$%
\textsl{\ be the associated Cartan matrix.} \textsl{For }$u,v\in W$\textsl{\
we have}
\end{quote}

\begin{center}
$(-1)^{l(w)}K_{u,v}^{w}=\Delta _{A_{w}}[(\sum b_{I}(u)y_{I})(\sum
b_{L}(v)y_{L})]_{(m)}$

$-\sum\limits_{\substack{ l(u)+l(v)\leq l(x)  \\ x<w}}b_{[1,\cdots
,m]}(x)K_{u,v}^{x}$\textsl{,}
\end{center}

\begin{quote}
\noindent \textsl{where }$I,L\subseteq \lbrack 1,\cdots ,m]$\textsl{, and
where the numbers }$b_{K}(x)$\textsl{, }$K\subseteq \lbrack 1,\cdots ,m]$%
\textsl{, }$x\in W$\textsl{, are given as that in Lemma 5.2.}
\end{quote}

\bigskip

Proceeding to the general case let $H\subset G$ be the centralizer of a
one--parameter subgroup in $G$. Take a maximal torus $T\subset H$ and
consider the standard fibration $p:G/T\rightarrow G/H$. It is well known
that (cf. [PR$_{2}$, Proposition 1.6])

\begin{quote}
\textsl{The induced ring map }$p^{!}:K(G/H)\rightarrow K(G/T)$\textsl{\ is
injective and satisfies }$p^{!}[\Omega _{w}(H)]=\Omega _{w}$\textsl{, }$w\in 
\overline{W}$\textsl{.}
\end{quote}

\noindent Consequently one gets

\begin{quote}
\textbf{Corollary 2.} $K_{u,v}^{w}(H)=K_{u,v}^{w}$\textsl{\ for }$u,v,w\in 
\overline{W}$\textsl{.}
\end{quote}

\begin{center}
\textbf{References}
\end{center}

\begin{enumerate}
\item[{[A]}] M. Atiyah, $K$-theory, Benjamin, Inc., New York-Amsterdam, 1967

\item[{[AH]}] M. Atiyah and F. Hirzebruch, Vector bundles and homogeneous
spaces. 1961 Proc. Sympos. Pure Math., Vol. III pp. 7--38 American
Mathematical Society, Providence, R.I.

\item[{[Br$_{1}$]}] M. Brion, Positivity in the Grothendieck group of complex
flag varieties. Special issue in celebration of Claudio Procesi's 60th
birthday. J. Algebra 258 (2002), no. 1, 137--159.

\item[{[Br$_{2}$]}] M. Brion, Lectures on the geometry of flag varieties, 
\textsl{preprint available on} arXiv: math.AG/0410240

\item[{[BS]}] R. Bott and H. Samelson, Application of the theory of Morse to
symmetric spaces, Amer. J. Math., Vol. LXXX, no. 4 (1958), 964-1029.

\item[{[Bu]}] A.S. Buch, A Littlewood-Richardson rule for the $K$-theory of
Grassmannians, Acta Math. 189 (2002), no. 1, 37--78.

\item[{[C]}] C. Chevalley, Les classes d'equivalence rationelle, I,II,
Seminaire C. Chevalley, Anneaux de Chow et applications (mimeographed
notes), Paris, 1958.

\item[{[Ch]}] C. Chevalley, Sur les D\'{e}compositions Cellulaires des
Espaces $G/B$, in Algebraic groups and their generalizations: Classical
methods, W. Haboush ed. Proc. Symp. in Pure Math. 56 (part 1) (1994), 1-26.

\item[{[D]}] M. Demazure, D\'{e}singularisation des vari\'{e}t\'{e}s de
Schubert g\'{e}n\'{e}ralis\'{e}es, Ann. Sci. \'{E}cole. Norm. Sup. (4)
7(1974), 53-88.

\item[{[Du$_{1}$]}] H. Duan, The degree of a Schubert variety, Adv. Math.,
180(2003), 112-133.

\item[{[Du$_{2}$]}] H. Duan, Multiplicative rule of Schubert classes, Invent.
Math.159(2005), 407-436.

\item[{[DZ$_{1}$]}] H. Duan and Xuezhi Zhao, A unified formula for Steenrod
operations in flag manifolds, \textsl{preprint available on} arXiv:
math.AT/0306250

\item[{[DZ$_{2}$]}] H. Duan and Xuezhi Zhao, Algorithm for multiplying
Schubert classes, \textsl{preprint available on} arXiv: math.AG/0309158.

\item[{[FL]}] W. Fulton and A. Lascoux, A Pieri formula in the Grothendieck
ring of a flag bundle, Duke Math. J. 76 (1994) 711--729.

\item[{[GR]}] S. Griffeth, and A. Ram, Affine Hecke algebras and the Schubert
calculus, European J. Combin. 25 (2004), no. 8, 1263--1283.

\item[{[HPT]}] W. Y. Hsiang, R. Palais and C. L. Terng, The topology of
isoparametric submanifolds, J. Diff. Geom., Vol. 27 (1988), 423-460.

\item[{[Hu]}] J. E. Humphreys, Introduction to Lie algebras and
representation theory, Graduated Texts in Math. 9, Springer-Verlag New York,
1972.

\item[{[KK]}] B. Kostant and S. Kumar, $T$-equivariant $K$-theory of
generalized flag varieties. J. Differential Geom. 32 (1990), no. 2, 549--603.

\item[{[L]}] C. Lenart, A $K$-theory version of Monk's formula and some
related multiplication formulas. J. Pure Appl. Algebra 179 (2003), no. 1-2,
137--158.

\item[{[LS]}] P. Littelmann and C.S. Seshadri, A Pieri-Chevalley type formula
for $K(G/B)$ and standard monomial theory. Studies in memory of Issai Schur
(Chevaleret/Rehovot, 2000), 155--176, Progr. Math., 210, Birkhauser Boston,
Boston, MA, 2003.

\item[{[M]}] O. Mathieu, Positivity of some intersections in $K_{0}(G/B)$, J.
Pure Appl. Algebra 152 (2000) 231--243.

\item[{[PR$_{1}$]}] H. Pittie and A. Ram, A Pieri-Chevalley formula in the $K$%
-theory of a $G/B$-bundle, Electron. Res. Announc. Amer. Math. Soc. 5 (1999)
102--107.

\item[{[PR$_{2}$]}] H. Pittie and A. Ram, A Pieri-Chevalley formula for $%
K(G/B)$, \textsl{preprint available on} arXiv: math.RT/0401332.

\item[{[W$_{1}$]}] M. Willems, Cohomologie et $K$--th\'{e}orie \'{e}%
quivariantes des tours de Bott et des vari\'{e}t\'{e}s de drapeaux.
Application au calcul de Schubert, \textsl{preprint available on} arXiv:
math.AG/0311079.

\item[{[W$_{2}$]}] M. Willems, K-th\'{e}orie \'{e}quivariante des vari\'{e}t%
\'{e}s de Bott-Samelson. Application \`{a} la structure multiplicative de la
K-th\'{e}orie \'{e}quivariante des vari\'{e}t\'{e}s de drapeaux, \textsl{%
preprint available on} arXiv: math.AG/0412152.
\end{enumerate}

\end{document}